\documentclass{svmult} 
\usepackage{mathptmx}
\usepackage{helvet}
\usepackage{courier}
\usepackage{graphicx}
\usepackage{makeidx}
\usepackage{multicol}
\usepackage{footmisc}
\usepackage{epsfig}
  \usepackage{t1enc}
      \usepackage[latin1]{inputenc}
      \usepackage[english]{babel}
          \usepackage{latexsym}
  \usepackage{amssymb}
 \begin{document}
\title*{SOME RESULTS AND PROBLEMS FOR ANISOTROPIC RANDOM WALKS ON THE PLANE}
\author{Endre Cs\'{a}ki\and Ant\'onia F\"oldes\and P\'al R\'ev\'esz}
\institute{Endre Cs\'aki \at Alfr\'ed R\'enyi Institute of Mathematics, 
Hungarian Academy of Sciences,\\ Budapest, P.O.B. 127, H-1364, Hungary,\\
\email{csaki.endre@renyi.mta.hu}
\and
Ant\'{o}nia F\"{o}ldes 
\at Department of Mathematics, College of Staten 
Island,\\ CUNY, 2800 Victory Blvd., Staten Island, New York 10314, U.S.A.\\ 
\email{antonia.foldes@csi.cuny.edu}
\and
P\'al R\'ev\'esz
\at Institut f\"ur Statistik und Wahrscheinlichkeitstheorie,\\ 
Technische
 Universit\"at Wien, Wiedner Hauptstrasse 8-10/107\\ A-1040 Vienna,
 Austria, \\
\email{revesz.paul@renyi.mta.hu}}

\renewcommand{\theequation}{\thesection.\arabic{equation}}
\newcommand{\dto}{\stackrel{{d}}{\rightarrow}}
\newtheorem{propo}{Proposition}[section]
\newtheorem{thm}{Theorem}[section]
\newtheorem{cjt}{Conjecture}[section]
 \def\begg{\begin{equation}}
 \def\endd{\end{equation}}
 \def\ep{\varepsilon}
 \def\noo{n\to\infty}
 \def\al{\alpha}
 \def\bE{\bf E}
 \def\be{\bf e}
 \def\bu{\bf u}
 \def\bv{\bf v}
 \def\bp{\bf P}
 \def\pe{\bf p}
 \def\bX{\bf X}
 \def\bY{\bf Y}
\def\pvu{p(\bf v,\bf u)}

\titlerunning{Anisotropic random walk}
\maketitle

\abstract{This is an expository paper on the asymptotic results concerning path
behaviour of the anisotropic random walk on the two-dimensional square
lattice $\mathbb{Z}^2$. In recent years Mikl\'os and the authors of the 
present  paper investigated the properties of this random walk concerning 
strong approximations, local times and range. We give a survey of these results
together with some further problems.}

\keywords{anisotropic random walk;
strong approximation; 2-dimensional Wiener process; local time; laws of
the iterated logarithm.} 
\vspace{.1cm}

\noindent
{\bf Mathematics Subject Classification (2010):}  
Primary 60F17, 60G50; Secondary 60F15, 60J10, 60J65

 \section{Introduction}
 \renewcommand{\thesection}{\arabic{section}} \setcounter{equation}{0}
 \setcounter{thm}{0} \setcounter{lemma}{0}

We consider random walks on the square lattice $\mathbb{Z}^2$ 
with possibly unequal symmetric horizontal and vertical step
probabilities, so that these probabilities can only depend on the value of
the vertical coordinate. In particular, if such a random walk is situated
at a site on the horizontal line $y=j\in \mathbb{Z}$, then at the next
step it moves with probability $p_j$ to either vertical neighbor, and with
probability $1/2-p_j$ to either horizontal neighbor.
More formally, consider the random walk
$\{{\bf C}(N)=\left(C_1(N),C_2(N)\right);\, N=0,1,2,\ldots\}$
on $\mathbb{Z}^2$ with the transition probabilities
\begin{eqnarray} \nonumber
{\bf P}({\bf C}(N+1)=(k+1,j)|{\bf C}(N)=(k,j))\\
\nonumber ={\bf P}
({\bf C}(N+1)=(k-1,j)|{\bf C}(N)=(k,j))=\frac{1}{2}-p_j,\\
\nonumber
{\bf P} ({\bf C}(N+1)=(k,j+1)|{\bf C}(N)=(k,j))\\
={\bf P}
({\bf C}(N+1)=(k,j-1)|{\bf C}(N)=(k,j))=p_j, \label{trans}
\end{eqnarray}
for $(k,j)\in{\mathbb Z^2}$, $N=0,1,2,\ldots$ We assume throughout the
paper that $0<p_j\leq 1/2$ and $\min_{j\in\mathbb Z} p_j<1/2$. Unless
otherwise stated we assume also that ${\mathbf C}(0)=(0,0)$.
This model has a number of physical applications and the topic has a broad
literature. We refer to Silver {\it et al.} \cite{SSL}, Seshadri {\it et al.} 
\cite{SLS},  Shuler
\cite{SH}, Westcott \cite{WE}, where certain properties of this random
walk were studied under various conditions. Heyde \cite{H} proved an
almost sure approximation for $C_2(\cdot)$ under the condition

\begg n^{-1}\sum_{j=1}^n p_j^{-1}=2\gamma+o(n^{-\eta}), \qquad
n^{-1}\sum_{j=1}^n p_{-j}^{-1}=2\gamma+o(n^{-\eta})\label{h} \endd
as $n\to\infty$ for some constants $\gamma$,  $1<\gamma< \infty$ and
$1/2<\eta<\infty .$

 Heyde {\it et al.} \cite{HWW} treated the case when conditions
similar to (\ref{h}) are assumed but $\gamma$ can be different for the two
parts of (\ref{h}) and obtained almost sure convergence to
the so-called oscillating Brownian motion. In Heyde \cite{H93} limiting
distributions were given for $\mathbf{C}(\cdot)$ under the condition
(\ref{h}) but without remainder. Den Hollander \cite {DH}
proved strong approximations for $\mathbf{C}(\cdot)$ in the case when
$p_j$-s are random variables with values $1/4$ and $1/2$. Roerdink and
Shuler \cite{RS} proved some asymptotic properties, including local limit 
theorems,
under certain conditions. For more detailed history see \cite{DH}.

First we give a general construction and discuss the issue of recurrence and 
transience of
this random walk. In Section 2 we discuss strong approximations of the random
walk $\{{\bf C}(N),\, N=0,1,\ldots\}$. Section 3 treats the local time and
in Section 4 some results on the range will be given.

\subsection{General construction}

Suppose that in a probability space we have
two independent simple symmetric random walks, i.e.,
$$
S_1(n),\, n=0,1,2,\ldots,\qquad S_2(n),\, n=0,1,2,\ldots,
$$
where $S_1(0)=S_2(0)=0$, $S_i(\cdot)$ are sums of i.i.d. random variables
each taking the values $1$ and $-1$ with probability $1/2$. The local times 
of $S_i$ are defined by
$$
\xi_i(j,n)=\#\{0\leq k\leq n: S_i(k)=j\}, \quad j\in \mathbb Z, \quad
n=0,1,2,\ldots 
$$

Moreover, on
the same probability space we have a double array of independent geometric
random variables
$$
G_i^{(j)},\, i=1,2,\ldots,\, j\in\mathbb Z
$$
with distribution
$$
\mathbf{P}(G_i^{(j)}=k)=2p_j(1-2p_j)^k,\,\, k=0,1,2,\ldots \label{geo}
$$

\bigskip
We now construct our walk $\mathbf{C}(N)$ as follows. We will take all
the horizontal steps consecutively from $S_1(\cdot)$ and all the vertical
steps consecutively from $S_2(\cdot).$ First we will take some horizontal
steps from $S_1(\cdot)$, then exactly one vertical step from $S_2(\cdot),$
then again some horizontal steps from $S_1(\cdot)$ and exactly one
vertical
step from $S_2(\cdot),$ and so on. Now we explain how to get the number of
horizontal steps on each occasion.
Consider our walk starting from the origin proceeding first horizontally
$G_1^{(0)}$ steps (note that $G_1^{(0)}=0$ is possible with probability
$2p_0$), after which it takes exactly one vertical step, arriving
either to the level $1$ or $-1$, where it takes $G_1^{(1)}$ or
$G_1^{(-1)}$ horizontal steps (which might be no steps at all) before
proceeding with another vertical step. If this step carries the walk to
the level $j$, then it will take $G_1^{(j)}$ horizontal steps, if this is
the first visit to level $j,$ otherwise it takes $G_2^{(j)}$ horizontal
steps. In general, if we finished the $k$ -th vertical step and arrived to
the level $j$ for the $i$-th time, then it will take $G_i^{(j)}$
horizontal steps.

In this paper $N$ will denote the number of steps of the walk
out of which $H_N$ denotes the number of horizontal steps and $V_N=n$
the number of vertical steps, i.e., $H_N+V_N=N$. Then clearly
$$
{\bf C}(N)=(C_1(N),C_2(N))=(S_1(H_N),S_2(V_N)).
$$

\subsection{Recurrence, transience}

Our result on recurrence is a simple application of the celebrated
Nash-Williams theorem \cite{NW}. To state this result we need some definitions.
Consider a Markov chain $(\mathbf{X},\mathbf{Y}, p)$ with countable state
space $\mathbf{X}$,  process $\mathbf{Y}$ and transition probabilities
$p(\mathbf{u},\mathbf{v}).$ The chain is reversible if there exist strictly
positive weights $\pi_{\mathbf{u}} $ for all $\mathbf{u} \in \bX$ such that

\begg \pi_{\mathbf{u}}p(\mathbf{u},\mathbf{v})=
\pi_{\mathbf{v}} p(\mathbf{v},\mathbf{u}).\endd
If the chain is reversible we will use
 $$a(\mathbf{u},\mathbf{v}):=\pi_{\mathbf{u}}p(\mathbf{u},\mathbf{v}).$$

Obviously the above defined anisotropic walk is a Markov chain on the state
space $\mathbf{X}=\mathbb{Z}^2$, with the transition probabilities defined in
(\ref{trans}). Furthermore, this Markov chain is reversible with the strictly
positive weights
$$\displaystyle{\pi(k,j)=\frac{1}{p_j}}$$
and
\begin{eqnarray} \nonumber
a\left((k,j),(k,j+1)\right)&=&a\left((k,j),(k ,j-1)\right)=1\\
a\left((k,j),(k+1,j)\right)&=&a\left((k,j),(k-1 ,j)\right)=
\frac{1}{2p_j}-1  \label{auv}
\end{eqnarray}
(and for non nearest neighbor sites $a(.,.) =0$).
This Markov chain is also time homogeneous, irreducible, i.e. it is possible
to get to any state from any state with positive probability. The invariant
measure is given by
\begin{equation}
\mu(k,j)=\pi(k,j)=\frac{1}{p_j}, \quad (k,j)\in \mathbb{Z}^2,
\label{invariant}
\end{equation}
i.e.,
$$
\mu(\mathbf{u})=\sum_{\mathbf{v}}\mu(\mathbf{v})p(\mathbf{v},\mathbf{u}),
$$
where the summation goes for the four nearest neighbors  of $\mathbf{u}$.

\smallskip\noindent
{\bf Theorem A} (Nash-Williams \cite{NW}): \textit{Suppose that
$(\mathbf{X},\mathbf{Y}, p)$ is a reversible Markov chain and that
${\bX}=\bigcup_{k=0}^{\infty}\Lambda^k$ where $\Lambda^k$ are disjoint.
Suppose further that  ${\bu} \in \Lambda^k$ and $a(\mathbf{u},\mathbf{v})>0$
together imply  that $\mathbf{v}
\in \Lambda^{k-1}\bigcup\Lambda^k\bigcup\Lambda^{k+1},$ and that for each
$k$ the sum $\displaystyle{\sum_{ \mathbf{u} \in \Lambda^k, \mathbf{v} \in
\mathbf{X}} a(\mathbf{u},\mathbf{v})< \infty}.$
Let $[\Lambda^{k},\Lambda^{k+1}]$ denote the set of pairs
$(\mathbf{u},\mathbf{v})$ such that $\mathbf{u} \in \Lambda^{k}$ and
$\mathbf{v} \in \Lambda^{k+1}.$  The Markov chain is recurrent if}
\begg
\sum_{k=0}^{\infty}\left( \sum_{({\mathbf{u},\mathbf{v})}
\in [\Lambda^{k},\,\Lambda^{k+1}]} a(\mathbf{u},\mathbf{v})\right)^{-1}
=\infty. \label{reccr}
\endd

To apply this theorem, let $\Lambda^k$ be  the set of $8k$ lattice points on
the square of width  $2k$, centered at the origin.  Furthermore, let
$[\Lambda^{k},\,\Lambda^{k+1}]$ be the set of $8k+4$ nearest neighbor pairs
(edges) between $\Lambda^{k}$ and $\Lambda^{k+1}.$

It is easy to see by (\ref{auv}) that the sum in (\ref{reccr}) is equal to
$$\sum _{k=0}^{\infty} \left(2\left(\sum_{j=-k}^k \left(\frac{1}{2p_j}-1\right)
+\sum_{j=-k}^k 1\right)\right)^{-1}=\sum _{k=0}^{\infty}\left(\sum_{j=-k}^k
\frac{1}{p_j}\right)^{-1}.$$
So we got the following result.
\begin{thm}
 The anisotropic walk is recurrent if
\begg
\sum _{k=0}^{\infty}\left(\sum_{j=-k}^k \frac{1}{p_j}\right)^{-1}=\infty.
\label{anirec}
\endd
\end{thm}
As a simple consequence, if $\min_{j\in \mathbb Z} p_j>0$, then the
anisotropic walk is recurrent.

It is an intriguing question whether the converse of this
statement is true as well. That is to say, is it true that

\begin{cjt} If
\begg
\sum _{k=0}^{\infty}\left(\sum_{j=-k}^k \frac{1}{p_j}\right)^{-1}<\infty,
\label{con}
\endd
then the anisotropic walk is transient.
\end{cjt}

We can't prove this conjecture, but a somewhat weaker result is true.

\begin{thm} Assume that
\begin{equation}
\sum_{j=-k}^k\frac{1}{p_j}=Ck^{1+A}+O(k^{1+A-\delta}),\quad k\to\infty
\label{transient}
\end{equation}
for some $C>0$, $A>0$ and $0<\delta\leq 1$. Then the anisotropic random walk is
transient.
\end{thm}

\noindent{\bf Proof.}
Consider the simple symmetric random walk $S_2(\cdot)$ of the vertical steps
and let $\xi_2(\cdot,\cdot)$ be its local time, and $\rho_2(\cdot)$ be its
return time to zero. Consider the anisotropic random walk of $N$ steps, where
$N=N(m)$ is the time of $m$-th return of $S_2(\cdot)$ to zero, i.e., let
$V_N={\rho_2(m)}$.  

First we give a lower bound for the number of the horizontal steps $H_N$.

\begin{lemma}
For small enough $\varepsilon>0$ we have almost surely for large enough $m$
\begin{equation}
H_N=H_{N(m)}\geq m^{1+(1-\varepsilon)(A+1)}.
\end{equation}
\end{lemma}

\noindent{\bf Proof.} For simplicity in the proof, we denote $\xi_2$ by
$\xi$ and $\rho_2$ by $\rho$.
From the construction in Section 1.1 it can be seen that the number of
horizontal steps up to the $m$-th return to zero by the vertical component
is given by
$$
H_N=\sum_{j=-\infty}^{\infty}\sum_{i=1}^{\xi(j,\rho(m))}G_i^{(j)},
$$
where $G_i^j$ are as in Section 1.1. Since $\rho(m)\geq m^{2-\varepsilon}$
for small $\varepsilon>0$ and large enough $m$ almost surely, it follows
from the stability of local time (see \cite{RE}, Theorem 11.22, p. 127),
that for any $\varepsilon>0$, $|j|\leq m^{1-\varepsilon}$ we have
$$
(1-\varepsilon)m\leq \xi(j,\rho(m))\leq (1+\varepsilon)m
$$
almost surely for large enough $m$. Hence
$$
H_N\geq\sum_{|j|\leq m^{1-\varepsilon}}
\sum_{i=1}^{(1-\varepsilon)m}G_i^{(j)}=:U_m.
$$
 We consider the expectation of $U_m$ and show that
the other terms are negligible. We have
$$
{\bf E}G_i^{(j)}=\frac{1-2p_j}{2p_j},
$$
$$
Var G_i^{(j)}=\frac{1-2p_j}{4p_j^2}.
$$
Hence by (\ref{transient}) we get
$$
{\bf E}U_m=m(1-\varepsilon)\sum_{|j|\leq m^{1-\varepsilon}}
\frac{1-2p_j}{2p_j}\geq cm^{1+(1-\varepsilon)(A+1)},
$$
where $c>0$  is a constant. In what follows  the value of such  a $c$
might change from line to line. We have
$$
Var U_m=m(1-\varepsilon)\sum_{|j|\leq m^{1-\varepsilon}}
\frac{1-2p_j}{4p_j^2}.
$$
It follows from (\ref{transient}) that $\frac{1}{2p_j}
\leq c|j|^{1+A-\delta}$,
hence
$$
Var U_m
\leq cm(1-\varepsilon)\sum_{|j|\leq m^{1-\varepsilon}}
\frac{|j|^{1+A-\delta}}{p_j}\leq cm^{1+(1-\varepsilon)(2+2A-\delta)}.
$$
By Chebyshev inequality
$$
{\bf P}(|U_m-{\bf E}U_m|\geq m^{(1-\varepsilon)(A+2)})\leq
c\frac{m^{3+2A-2\varepsilon(1+A)-(1-\varepsilon)\delta}}
{m^{2(1-\varepsilon)(A+2)}}=cm^{-1-(1-\varepsilon)\delta+2\varepsilon}
$$
which, by choosing $\varepsilon$ small enough, is summable. Hence, as
$m\to\infty$,
$$
U_m={\bf E}U_m+O(m^{(1-\varepsilon)(A+2)})\quad a.s.,
$$
consequently
$$
H_N\geq U_m\geq cm^{1+(1-\varepsilon)(1+A)}
$$
almost surely for large $m$. $\Box$

\begin{lemma}
Let $S(\cdot)$ be a simple symmetric random walk and let $r(m)$ be a
sequence of integer valued random variables independent of $S(\cdot)$ and 
such that $r(m)\geq m^{\beta}$ almost surely for large $m$ with $\beta>2$. 
Then with small enough $\varepsilon>0$ we have
$$
|S(r(m))|\geq m^{\beta/2-1-\varepsilon}
$$
almost surely for large $m$.
\end{lemma}
{\bf Proof.} 
From the local central limit theorem
$$
{\bf P}(S(k)=j)\leq \frac{c}{\sqrt{k}}
$$
for all $k\geq 0$ and $j\in \mathbb{Z}$ with an absolute constant $c>0$.
Hence
$$
{\bf P}\left(\frac{|S(k)|}{\sqrt{k}}\leq x\right)=
\sum_{|j|\leq x\sqrt{k}}P(S(k)=j)\leq cx,
$$
This remains true if $k$ is replaced by a random variable, independent of
$S(\cdot)$, e.g. $k=r(m)$, i.e. we have
$$
{\bf P}\left(\frac{|S(r(m))|}{\sqrt{r(m)}}\leq
\frac{1}{m^{1+\varepsilon}}\right)
\leq c\frac{1}{m^{1+\varepsilon}},
$$
consequently by Borel-Cantelli Lemma
$$
|S(r(m))|\geq \frac{\sqrt{r(m)}}{m^{1+\varepsilon}}\geq
m^{\beta/2-1-\varepsilon},
$$
almost surely for all large enough $m$. This completes the proof of the 
Lemma. $\Box$

Applying the two lemmas with $r(m)=H_{N(m)}$, we get
$$
|S_1(H_N)|\geq m^{A/2-\varepsilon A/2-3\varepsilon/2}=m^{\gamma}
$$
with $\gamma>0$ by choosing $\varepsilon>0$ small enough.
It follows that for large $N$, $S_1(H_N)$ almost surely can't be equal to zero.

Let
$$
A_m:= \{  \exists j, \rho_2(m)<j <\rho_2(m+1)\quad
{\rm such\,\, that}\quad \mathbf{C}(j)=(0,0) \}.
$$
Clearly  $A_m$ could only happen if from $\rho_2(m)$ to $\rho_2(m+1)$ the
walk only steps horizontally  (if it makes one vertical step the return to 
the origin  could only
happen after or at $\rho_2(m+1)$). Thus by our lemmas in order that $A_m$ 
could happen, the walk needs to have  at least
$m^{\gamma}$ consecutive steps on the $x$-axis, thus
$$\sum_{m}^{\infty} \mathbf{P}(A_m)\leq
\sum_{m}^{\infty}(1/2-p_0)^{m^{\gamma}}<\infty.$$
So the anisotropic random walk cannot return to zero infinitely often with
probability 1, it is transient. This proves the theorem. $\Box$

\section{Strong approximations}
 \renewcommand{\thesection}{\arabic{section}} \setcounter{equation}{0}

In this section we present results concerning strong approximations of the
two-dimensional aniso\-tropic random walks. Of course, the results are
different in the various cases, and in some cases the problem is open.
We also mention weak convergence results available in the literature.
First we describe the general method how to obtain these strong approximations.

Assume that our anisotropic random walk is constructed on a probability
space as described in Section 1.1, and in accordance with Theorems 6.1 and 
10.1 of R\'ev\'esz \cite{RE} we may assume that on the same probability space
there are also two independent standard Wiener processes (Brownian motions)
$W_i(\cdot),\, i=1,2$ with local times $\eta_i(\cdot,\cdot)$ such that
for all $\varepsilon>0$, as $n\to\infty$, we have
$$
S_i(n)=W_i(n)+O(n^{1/4+\varepsilon}) \quad a.s.
$$
and
$$
\xi_i(j,n)=\eta_i(j,n)+O(n^{1/4+\varepsilon}) \quad a.s.
$$
Then
$$
C_1(N)=S_1(H_N)=W_1(H_N)+O(H_N^{1/4+\varepsilon}) \quad a.s.,
$$
and
$$
C_2(N)=S_2(V_N)=W_2(V_N)+O(V_N^{1/4+\varepsilon}) \quad a.s.,
$$
if $H_N\to\infty$ and $V_N\to\infty$ as $N\to\infty$, almost surely.

So we have to give reasonable approximations to $H_N$ and $V_N$, or at
least to one of them, and use $V_N+H_N=N$.

It turned out that in many cases treated, the following is a good
approximation of $H_N$.
$$
H_N\sim \sum_j\sum_{i=1}^{\xi_2(j,n)}G_i^{(j)}
\sim \sum_j\xi_2(j,n)\frac{1-2p_j}{2p_j},
$$
with $n=V_N$. $H_N$ and the double sum above are not necessarily equal, since 
the last geometric variable might be truncated in $H_N$.
So we have to investigate the additive functional

$$
A(n)=\sum_j f(j)\xi_2(j,n)=\sum_{k=0}^nf(S_2(k)),\quad
f(j)=\frac{1-2p_j}{2p_j}
$$
of the vertical component and approximate it by the additive functional
of $W_2(\cdot)$

$$
B(t)=\int_{-\infty}^{\infty}f(x)\eta_2(x,t)\, dx=
\int_0^tf(W_2(s))\, ds,
$$
where between integers we define $f(x)=f(j)$, $j\leq x<j+1$.

In certain cases the following Lemma, equivalent to Lemma 2.3 of Horv\'ath 
\cite{HO}, giving strong approximation of additive functionals, may be useful. 

\noindent
{\bf Lemma A} \textit{Let $g(t)$ be a non-negative function such that
$g(t)=g(j)$, $j\leq t<j+1$, for $j\in\mathbb{Z}$ and assume that
$g(t)\leq c(|t|^\beta+1)$ for some $0<c$ and $\beta\geq 0$. Then
$$
|\sum_{j=0}^n g(S_2(j))-\int_0^n g(W_2(s)) ds|=
o(n^{\beta/2+3/4+\varepsilon})\quad a.s.
$$
as} $n\to\infty$.

Now let us introduce the notations
$$
\sum_{j=1}^kf(j)=b_k,\quad \sum_{j=1}^kf(-j)=c_k
$$

The next assumption is a reasonable one used in the literature: as $k\to\infty$
\begin{equation}
b_k=(\gamma-1) k^\alpha +o(k^{\alpha-\delta})
\label{ass1}
\end{equation}
\begin{equation}
c_k=(\gamma-1) k^\alpha +o(k^{\alpha-\delta})
\label{ass2}
\end{equation}
with some $\gamma\geq 1$, $\alpha\geq 0$ and $\delta>0$. Observe that  
(\ref{h}) is a particular case with $\alpha=1$.

We consider the following cases:

(1) $\alpha=0$

(2) $0<\alpha<1$

(3) $\alpha=1$

(4) $\alpha>1$

(5) nonsymmetric case, i.e. the constants $\gamma$ in (\ref{ass1}) and
(\ref{ass2}) are different.

\subsection{The case $\alpha=0$}

The most interesting and well established case is the so-called comb
structure, i.e., $p_0=1/4$, $p_j=1/2,\, \, j=\pm 1,\pm 2,\ldots.$
It follows from Theorem 1.1 that the random walk in this case is recurrent. 
We note in passing the interesting result of Krishnapur-Peres \cite{KP}: 
two independent random walks on the comb meet only finitely often with 
probability 1.

For random walk on comb we refer to Weiss and Havlin \cite{WH}, Bertacchi
and Zucca \cite{BZ} and references given there.
The following result on weak convergence was established by Bertacchi
\cite{BE}.

\medskip\noindent{\bf Theorem B} 
$$
\left(\frac{C_1(Nt)}{N^{1/4}},\frac{C_2(Nt)}{N^{1/2}};\, t\geq 0\right)
{\buildrel{\rm Law}\over\longrightarrow}\, (W_1(\eta_2(0,t)), W_2(t);\,
t\geq 0), \quad N\to\infty.
$$

Strong approximation was given in Cs\'aki {\it et al.} \cite{CCFR09}.
\begin{thm}
On an appropriate probability space we have
$$
N^{-1/4}|C_1(N)-W_1(\eta_2(0,N))|+N^{-1/2}|C_2(N)-W_2(N)|
=O(N^{-1/8+\varepsilon})\quad {a.s.},
$$
as $N\to\infty$.
\end{thm}

We have the following consequences.
\begin{corollary}
$$
\limsup_{N\to\infty}\frac{C_1(N)}{2^{5/4}3^{-3/4}N^{1/4}(\log\log N)^{3/4}}
=1\quad a.s.
$$
$$
\limsup_{N\to\infty}\frac{C_2(N)}{(2N\log\log N)^{1/2}}=1\quad a.s.
$$
\end{corollary}

For general results in the case $\alpha=0$ we just remark
that in this case $\bar f=\sum_jf(j)$ is convergent, then by our assumptions
its terms are non-negative and at least one of them is strictly positive,
hence $\bar f>0$. By the ratio ergodic theorem (cf., e.g., Theorem 3.6 
in Revuz \cite{REVUZ}) 
 
$$
A(n)\sim \bar f\xi_2(0,n), \quad \bar f=\sum_j f(j)=2(\gamma-1)+f(0),
$$
almost surely, as $n\to\infty$, hence
$$
A(n)=O((n\log\log n)^{1/2}) \quad {\rm a.s.},\, n\to\infty.
$$
Let
$$H_N^+=\sum_j\sum_{i=1}^{\xi_2(j,n)}G_i^{(j)},
\quad H_N^-=\sum_j\sum_{i=1}^{\xi_2(j,n)-1}G_i^{(j)}.
$$
Obviously, $H_N^-\leq H_N\leq H_N^+$. Consider $H_N^+$, which is a (random)
sum of independent random variables. Under the condition
${\cal F}=\{S_2(k),\, k\geq 0\}$ we have
$$
E(H_N^+|{\cal F})=\sum_jf(j)\xi_2(j,n)=A(n)
$$
$$
Var(H_N^+|{\cal F})=\sum_j\frac{f(j)}{2p_j}\xi_2(j,n).
$$
It is easy to see that the sum $\sum_jf(j)/2p_j$ is also convergent, hence
$$
Var(H_N^+|{\cal F})\sim c\xi(0,n)
$$
with some $c>0$.
Now apply Theorem 6.17 in Petrov \cite{PE} saying that for sums of 
independent (not necessary identically distributed) random variables we have
$$
\sum_i X_i=\sum_iEX_i+
o\left(\left(\sum_iVar X_i\right)^{1/2+\varepsilon}\right)
$$
almost surely. Thus
$$
H_N^+=\bar f\xi_2(0,n)(1+o(1))=\bar f\xi_2(0,V_N)(1+o(1))
$$
almost surely as $N\to\infty$. Similar results are true for $H_N^-$, hence
also for $H_N$, i.e.
$$
H_N=\bar f\xi_2(0,n)(1+o(1))=\bar f\xi_2(0,V_N)(1+o(1)).
$$

Since $C_1(N)=S_1(H_N)$, using that $H_N=O((N\log\log N)^{1/2})$ and the
strong approximations of $S_1(\cdot),S_2(\cdot)$ by $W_1(\cdot),W_2(\cdot)$
and $\xi_2(0,\cdot)$ by $\eta_2(\cdot)$, we can obtain the following limit
distribution: as $N\to\infty$,
$$
\left(\frac{C_1(N)}{N^{1/4}}, \frac{C_2(N)}{N^{1/2}}\right)
\dto \left(W_1(\bar f\eta_2(0,1)),W_2(1)\right).
$$

Further results, like strong approximations, remain to be established in this 
case. 

\subsection{The case $0<\alpha<1$}
\bigskip
This is also a recurrent case, but approximations, limit theorems, etc. 
remain to be worked out in detail. We just note that 
from the law of the iterated logarithm for the local time we have
$$
A(n)=\sum_jf(j)\xi_2(j,n)\leq c(n\log\log n)^{(1+\alpha)/2},
$$
a.s., $n\to\infty$, hence the vertical part dominates, i.e., as
$N\to\infty$ we should have
$$
H_N=O((N\log\log N)^{(1+\alpha)/2})<<N \quad a.s., 
$$
and we expect that
$$
C_1(N)=W_1(Z(N))+O(N^{(1+\alpha)/4+\varepsilon})\quad a.s.,
$$
where
$$
Z(N)=\int_0^N f(W_2(s))\, ds,
$$
and for the vertical component
$$
C_2(N)=W_2(N)+O(N^{1/2-\varepsilon}), \quad a.s.
$$
as $N\to\infty$.

\subsection{The case $\alpha=1$}

Assume also that $\delta>1/2$, $\gamma>1$.

It can be seen from Thoerem 1.1 that the anisotropic random walk in this 
case is recurrent. 

\bigskip
Heyde \cite{H} gave the following strong approximation:

\medskip\noindent
{\bf Theorem C} {\it On an appropriate probability space we have
$$
\gamma^{1/2}C_2(N)=W_2(N(1+\varepsilon_N))+
O(N^{1/4}(\log N)^{1/2}(\log\log N)^{1/2})\quad a.s.
$$
as $N\to\infty$, where $W(\cdot)$ is a standard Wiener process
and $\lim_{N\to\infty}\varepsilon_N=0$ a.s.}

In another paper Heyde \cite{H93} gave weak convergence result for both
coordinates.

\medskip\noindent
{\bf Theorem D} 
$$
\left(\frac{C_1(N)}{N^{1/2}},\frac{C_2(N)}{N^{1/2}}\right)
\dto \left(W_1(1-\gamma^{-1}), W_2(\gamma^{-1})\right).
$$

Strong approximation result for both coordinates was given
in Cs\'aki {\it et al.} \cite{CCFR12a}:
\begin{thm} On an appropriate probability space we have for any $\varepsilon>0$
$$
\left|C_1(N)-W_1\left(\frac{(\gamma-1)N}{\gamma}\right)\right|+
\left|C_2(N)-W_2\left(\frac{N}{\gamma}\right)\right|
=O(N^{5/8-\delta/4+\varepsilon})\quad {a.s.},
$$
as $N\to\infty$.

Moreover, in the periodic case, when $p_j=p_{j+L}$\, for each
$j\in \mathbb{Z}$ and a fixed integer $L\geq 1$ we have
$$
\left|C_1(N)-W_1\left(\frac{(\gamma-1)N}{\gamma}\right)\right|+
\left|C_2(N)-W_2\left(\frac{N}{\gamma}\right)\right|
=O(N^{1/4+\varepsilon}) \quad {a.s.},
$$
as $N\to\infty$, where
$$
\gamma=\frac{\sum_{j=0}^{L-1}p_j^{-1}}{2L}.
$$
\end{thm}

Some consequences are the following laws of the iterated logarithm.
$$
\limsup_{N\to\infty}\frac{C_1(N)}{(N\log\log
N)^{1/2}}=\left(\frac{2(\gamma-1)}{\gamma}\right)^{1/2} \quad {a.s.}
$$
$$
\limsup_{N\to\infty}\frac{C_2(N)}{(N\log\log
N)^{1/2}}=\left(\frac{2}{\gamma}\right)^{1/2} \quad {a.s.}
$$

\subsection{The case $\alpha>1$}

In this case the random walk is transient by Theorem 1.2.

It is an open problem to give strong approximations in this case. Horv\'ath
\cite{HO} established weak convergence of $C_2(\cdot)$ to some time changed
Wiener process. We mention a particular case of his results, valid for all 
$\alpha>1$. 

Let
$$
I_\alpha(t)=\int_0^t |W_2(s)|^{\alpha-1}\, ds.
$$
$I_\alpha$ is strictly increasing, so we can define its inverse,
denoted by $J_\alpha$. Then we have

$$
\frac{C_2(Nt)}{N^{1/(1+\alpha)}}
{\buildrel{\rm Law}\over\longrightarrow}\, c_0W_2(J_\alpha(t))
$$
with some constant $c_0$. 

In this case the number of horizontal steps dominates the number of vertical 
steps, therefore $C_1(N)$ might be approximated by $W_1(N)$.

\subsection{Unsymmetric case}

Some weak convergence in this case was treated in Heyde {\it et al.} \cite{HWW}
and Horv\'ath \cite{HO}. Strong approximation in a particular case, the
so-called half-plane half-comb structure (HPHC) was given in Cs\'aki {\it et 
al.} \cite{CCFR12b}.

Let $p_j=1/4$, $j=0,1,2,\ldots$ and
$p_j=1/2$, $j=-1,-2,\ldots$, i.e., a square lattice on the upper
half-plane, and a comb structure on the lower half plane. Let
furthermore
$$
\alpha_2(t):=\int_0^t I\{W_2(s)\geq 0\}\, ds,
$$
i.e., the time spent by $W_2$ on the non-negative side during the interval
$[0,t]$. The process $\gamma_2(t):=\alpha_2(t)+t$ is strictly increasing,
hence we can define its inverse: $\beta_2(t):=(\gamma_2(t))^{-1}$.
\begin{thm}
On an appropriate probability space we have
$$
|C_1(N)-W_1(N-\beta_2(N))|+|C_2(N)-W_2(\beta_2(N))|
=O(N^{3/8+\varepsilon})\quad {a.s.},
$$
as $N\to\infty$.
\end{thm}

The following laws of the iterated logarithm hold:
\begin{corollary}
$$
\limsup_{t\to\infty}\frac{W_1(t-\beta_2(t))}{\sqrt{t\log\log t}}=
\limsup_{N\to\infty}\frac{C_1(N)}{\sqrt{N\log\log N}}=1
\quad a.s.,
$$
$$
\liminf_{t\to\infty}\frac{W_1(t-\beta_2(t))}{\sqrt{t\log\log t}}=
\liminf_{N\to\infty}\frac{C_1(N)}{\sqrt{N\log\log N}}=-1
\quad a.s.,
$$
$$
\limsup_{t\to\infty}\frac{W_2(\beta_2(t))}{\sqrt{t\log\log t}}=
\limsup_{N\to\infty}\frac{C_2(N)}{\sqrt{N\log\log N}}=1
\quad a.s.,
$$
$$
\liminf_{t\to \infty}\frac{W_2(\beta_2(t))}{\sqrt{t\log\log t}}=
\liminf_{N\to\infty}\frac{C_2(N)}{\sqrt{N\log\log N}}=-\sqrt{2}
\quad a.s.
$$
\end{corollary}

\section{Local time}
 \renewcommand{\thesection}{\arabic{section}} \setcounter{equation}{0}
 \setcounter{thm}{0} \setcounter{lemma}{0}

We don't know any general result about the
local time of the anisotropic walk. It would require to determine asymptotic
results or at least good estimations for the return probabilities, i.e.,
we would need local limit theorems. In fact,
we know such results in two cases: the periodic and the comb structure case.

\subsection{Periodic anisotropic walk}

In case of the periodic anisotropic walk, i.e., when $p_j=p_{j+L}$, for some
fixed integer $L\geq 1$ and $j=0,\pm 1,\pm 2,\ldots$ we know the following 
local limit theorem for the random walk denoted by ${\bf C^P}(\cdot)$.

\begin{lemma} As $N\to\infty$, we have
\begg
\mathbf{P}(\mathbf{C^P}(2N)=(0,0))\sim \frac{1}{4\pi N
p_0\sqrt{\gamma-1}}
\label{return}
\endd
with $\gamma=\sum_{j=0}^{L-1}p_j^{-1}/(2L)$.
\end{lemma}

The proof of this lemma is based on the work  of Roerdink and Shuler \cite{RS}.
It follows from this lemma, that the truncated Green function $g(\cdot)$ is
given by
$$
g(N)=\sum_{k=0}^N \mathbf{P}(\mathbf{C^P}(k)=(0,0))\sim \frac{\log
N}{4p_0\pi\sqrt{\gamma-1}},\qquad N\to\infty,
$$
which implies that our anisotropic random walk in this case is recurrent
and also Harris recurrent.

First, we define the local time by
\begg
\Xi((k,j),N)=\sum_{r=1}^N I\{\mathbf{C^P}(r)=(k,j)\},\quad
(k,j)\in{\mathbb Z}^2.
\label{loctime}
\endd
In the case when the random walk is (Harris) recurrent, then we have (cf.
e.g. Chen \cite{CX})
$$
\lim_{N\to\infty}\frac{\Xi((k_1,j_1),N)}{\Xi((k_2,j_2),N)}=
\frac{\mu(k_1,j_1)}{\mu(k_2,j_2)}\quad {a.s.},
$$
where $\mu(\cdot)$ is an invariant measure.
Hence by (\ref{invariant})
$$
\lim_{N\to\infty}\frac{\Xi((0,0),N)}{\Xi((k,j),N)}=\frac{p_j}{p_0}
\quad {a.s.}
$$
for $(k,j)\in\mathbb{Z}^2$ fixed.

Thus, using now $g(N)$, it follows from Darling and Kac \cite{DK} that we
have
\begin{corollary}
$$
\lim_{N\to\infty}\mathbf{P}\left(\frac{\Xi((0,0),N)}{g(N)}\geq x\right)=
\lim_{N\to\infty}
\mathbf{P}\left(\frac{4p_0\pi\sqrt{\gamma-1}\, \Xi((0,0),N)}{\log N}\geq
x\right)=e^{-x}
$$
for $x\geq 0$.
\end{corollary}

For a limsup result, via Chen \cite{CX} we conclude
\begin{corollary}
$$
\limsup_{N\to\infty}\frac{\Xi((0,0),N)}{\log N\log\log\log N}
=\frac{1}{4p_0\pi\sqrt{\gamma-1}}\quad{a.s.}
$$
\end{corollary}

\subsection{Comb}
Now we consider the case of the two-dimensional comb structure ${\mathbb C}^2$,
i.e., when $p_0=1/4$ and $p_j=1/2$ for $j=\pm 1,\pm 2,\ldots$

First we give the return probability from Woess \cite{WO}, p. 197:
$$
{\bf P}({\bf C}(2N)=(0,0))\sim\frac{2^{1/2}}{\Gamma(1/4)N^{3/4}},
\quad N\to\infty.
$$

This result indicates that the local time tipically is of order $N^{1/4}$.
In Cs\'aki {\it et al.} \cite{CCFR10} and \cite{CCFR11} we have shown the 
following results.

\begin{thm}
The limiting distribution of the local time is given by
$$
\lim_{N\to\infty}{\bf P}(\Xi((0,0),N)/N^{1/4}<x)=
{\bf P}(2\eta_1(0,\eta_2(0,1))<x)={\bf P}(2|U|\sqrt{|V|}<x),
$$
where $U$ and $V$ are two independent standard normal random variables.
\end{thm}

Concerning strong approximation, in Cs\'aki {\it et al.} \cite{CCFR11} we 
proved the following results.

\begin{thm}
On an appropriate probability space for the random
walk  \newline $\{{\bf C}(N)=(C_1(N),C_2(N)); N=0,1,2,\ldots\}$ on
${\mathbb C}^2$, one can construct two independent standard Wiener
processes $\{W_1(t);\, t\geq 0\}$, $\{W_2(t);\, t\geq 0\}$ with their
corresponding local time processes $\eta_1(\cdot,\cdot),
\eta_2(\cdot,\cdot)$ such that, as $N\to\infty$, we have for any $\delta>0$
$$
\sup_{x\in{\mathbb Z}} \left|\Xi((x,0),N)-2
\eta_1\left(x,\eta_2(0,N)\right) \right|=O(N^{1/8+\delta})\quad {a.s.}
$$
\end{thm}

The next result shows that on the backbone up  to $|x|\leq N^{1/4-\epsilon}$
we have uniformity, all the sites have approximately the same local time as
the origin. Furthermore if we consider a site on a tooth of the comb its
local time is roughly half of the local time of the origin. This is pretty
natural, as it turns out from the proof that on the backbone the number of
horizontal and vertical visits to any particular site are approximately equal.

\begin{thm}
On the probability space of {\rm Theorem 3.2}, as
$N\to \infty,$  we have for any $0<\varepsilon <1/4$
$$
\max_{|x|\leq N^{1/4-\varepsilon}}|\Xi((x,0),N)-\Xi((0,0),N)|=
O(N^{1/4-\delta}) \quad a.s.
$$
and
$$
\max_{0<|y|\leq N^{1/4-\varepsilon}}\max_{|x|\leq
N^{1/4-\varepsilon}}
|\Xi((x,y),N)-\frac{1}{2}\Xi((0,0),N)| =O(N^{1/4-\delta})\quad {a.s.},
$$
for any $0<\delta<\varepsilon/2$, where the maximum
is taken on the integers.
\end{thm}

It would be an interesting problem to investigate the local time for
$|y|>N^{1/4}$. We believe e.g. that the maximal local time taken for all 
$(x,y)\in \mathbb{Z}^2$ is of order $N^{1/2}$. Such results however remain 
to be established. 

One of our old results \cite{CFR97} describes the Strassen class of
$\eta_1(0,\eta_2(0,zt))$ as follows. This, combined with Theorems 3.2 and
3.3, allows us to conclude the corresponding Strassen class result for the
local times of the walk.
\begin{thm} The net
$$
\left\{\frac{\eta_1(0,\eta_2(0,zt))}{2^{5/4}3^{-3/4}t^{1/4}(\log\log
t)^{3/4}};\, 0\leq z\leq 1\right\}_{t\geq 3},
$$
as $t\to\infty$, is almost surely relatively compact in the space
$C([0,1],{\mathbb R})$ of continuous functions from $[0,1]$ to ${\mathbb
R},$ and the set of its limit points is the class of nondecreasing
absolutely continuous functions {\rm (}with respect to the Lebesgue
measure{\rm )} on $[0,1]$ for which
$$
 {\cal S}^*:\left\{ f(0)=0\,\, {\rm and}\,\,
\int_0^1|\dot{f}(x)|^{4/3}\, dx\leq 1\right\}.$$
\end{thm}

Some obvious consequences of these results are the following

\begin{itemize}
\item
$\displaystyle{\limsup_{t\to\infty}\frac{\eta_1(0,\eta_2(0,t))}{
t^{1/4}(\log\log t)^{3/4}}=\frac{2^{5/4}}{3^{3/4}}\quad {\rm a.s.}}$

 \item $\displaystyle{\limsup_{N\to\infty}\frac{\Xi((x,0),N)}{N^{1/4}(\log\log
N)^{3/4}}=\frac{2^{9/4}}{3^{3/4}}\quad {\rm a.s.}},$

 \item $ \displaystyle{\limsup_{N\to\infty}\frac{\Xi((x,y),N)}{N^{1/4}(\log\log
N)^{3/4}}=\frac{2^{5/4}}{3^{3/4}}
\quad {\rm a.s.}\,\, y\neq 0}.$
\end{itemize}

A beautiful classical result of
L\'evy, P.  \cite{L} reads as follows

\medskip\noindent
{\bf Theorem E } {\it Let $W(\cdot)$ be a standard Wiener
process with local time process $\eta(\cdot,\cdot)$.  The following
equality in distribution holds}:
$$
\{\eta(0,t),\, t\geq 0\}{\buildrel{d}\over =}
\{\sup_{0\leq s\leq t}W(s), \, t\geq 0\}.
$$
Consequently using a Hirsch type result of
Bertoin  \cite{BER},  we get

\begin{corollary} Let $\beta(t)>0,\, t\geq 0$, be a non-increasing
function. Then we have almost surely that
$$
\liminf_{t\to\infty}\frac{\eta_1(0,\eta_2(0,t))}{t^{1/4}\beta(t)}=0\quad
or\quad \infty
$$
according as the integral $\int_1^\infty \beta(t)/t\, dt$
diverges or converges.
\end{corollary}

So we also have

\begin{corollary}
 Let $\beta(n), n=1,2,\ldots$ be a non-increasing
sequence of positive numbers. Then, for any fixed
$(x,y) \in {\mathbb Z}^2,$   we have almost surely that
$$
\liminf_{n\to\infty}\frac{\Xi((x,y),n)}{n^{1/4}\beta(n)}=0\quad
or\quad \infty
$$
according as the series $\sum_1^\infty \beta(n)/n$ diverges or
converges.
\end{corollary}

Now we also might consider the behavior of the supremum of the local time 
over the backbone.
To this end we first  had to  prove the following pair of integral tests for 
the $\sup_{x\in {\mathbb R}}\eta_1(x,\eta_2(0,t))$ process.

\begin{thm}
 Let $f(t)>0$ be a non-decreasing function and put
$$
I(f):=\int_1^\infty \frac{f^2(t)}t
\exp\left(-\frac{3}{2^{5/3}}f^{4/3}(t)\right)\, dt.
$$
Then, as $t\to \infty$,
$$
{\bf P}(\sup_{x\in {\mathbb R}}\eta_1(x,\eta_2(0,t))>t^{1/4} f(t)\, \, i.o.)
=0\, \, \, or\, \, \, 1
$$
according as $I(f)$ converges or diverges.
\end{thm}
\begin{thm}
 Let $g(t)>0$ be a non-increasing function and
$$
J(g):=\int_1^\infty \frac{g^2(t)}{t}\, dt.
$$
Then, as $t\to \infty$,
$$
{\bf P}(\sup_{x\in {\mathbb R}}\eta_1(x,\eta_2(0,t))<t^{1/4}g(t)\, \, i.o.)
=0\, \, \, or\, \, \, 1
$$
according as whether $J(g)$ converges or diverges.
\end{thm}
The above theorems imply the following integral tests for 
$\sup_{x\in {\mathbb Z}}\Xi((x,0),n);$

\begin{thm}  Let $a(n)$ be a non-decreasing sequence. Then, as $n\to\infty$,
$$
{\bf P}(\sup_{x\in {\mathbb Z}}\Xi((x,0),n)>2n^{1/4}a(n)\,\, 
i.o.)=0\, \, \, or\, \, \, 1
$$
according as
$$
\sum_{n=1}^\infty \frac{a^2(n)}{n}\exp\left(-\frac{3a^{4/3}(n)}
{2^{5/3}}\right)<\infty\, \, \, or\, \, \, =\infty.
$$
\end{thm}
\begin{thm} Let $b(n)$ be a non-increasing sequence. Then,
as $n\to\infty$,

$$
{\bf P}(\sup_{x\in {\mathbb Z}}\Xi((x,0),n)<n^{1/4}b(n)\, \, \, i.o.)
=0\, \, \, or\, \, \, 1
$$
according as
$$
\sum_{n=1}^\infty \frac{b^2(n)}{n}<\infty\, \, \, or\, \, \, =\infty.
$$
\end{thm}

 \section{Range}
\renewcommand{\thesection}{\arabic{section}} \setcounter{equation}{0}
\setcounter{thm}{0} \setcounter{lemma}{0}

The range of the anisotropic walk is defined in the usual way as

$$R(N)=\sum_{(k,j)\in \mathbb{Z}^2} I(\Xi((k,j),N)>0)$$
i.e., the number of distinct sites visited by the random walk up to time
$N,$ where $\Xi((k,j),N)$ is the local time of the  point $(k,j)$ at time $N.$

 We are not aware of any all embracing result about the range of the
anisotropic walk in general.
However the case of the periodic walk is completely understood.

Recall that the walk is periodic if
 $p_j=p_{j+L}$ for each $j\in \mathbb{Z},$
where $L\geq 1$ is a positive integer.
In this case it  is easy to see that
$$\gamma=\frac{\sum_{j=0}^{L-1} p_j^{-1}}{2L}.$$

 Roerdink and Shuler \cite{RS} gives  the asymptotic expected value of the
range:
$$
\mathbf{E}(R(N))\sim \frac{2\pi\sqrt{\gamma-1}}{\gamma}\frac{N}{\log N},
\quad N\to\infty.
$$

Moreover, it can be seen that our walk in this case is equivalent to the
so-called random walk with internal states, consequently, a law of large
numbers follows from N\'andori \cite{NA}

$$\lim_{N\to\infty}\frac{R(N)}{\mathbf{E}(R(N))}=\lim_{N\to\infty}
\frac{\gamma\, R(N)\log N}{2\pi\sqrt{\gamma-1}\, N}=1\qquad {a.s.}
$$

As a special case from these results we recover the well-known
Dvoretzky-Erd\H os \cite{DE} results for the simple random walk on the plane 
(without the remainder term), as for the plane $L=1$  and $\gamma=2$.
Thus we get 
$$
\mathbf{E}(R(N))\sim \frac{\pi N}{\log N},
\quad N\to\infty.
$$
and

$$\lim_{N\to\infty}\frac{R(N)}{\mathbf{E}(R(N))}=\lim_{N\to\infty}
\frac{ R(N)\log N}{\pi
 N}=1\qquad {a.s.}
$$

\noindent
{\bf Acknowledgement.} The authors thank the referees for valuable comments 
and suggestions. Research supported by PSC CUNY Grant, No. 68080-0043 and by
the Hungarian National Foundation for Scientific Research, No. K108615.

\end{document}